\documentclass[jcp,floats,amsmath,amssymb,onecolumn, floatfix]{revtex4-1}
\pdfoutput=1 

\usepackage{graphicx}
\usepackage{dcolumn}
\usepackage{bm}
\usepackage{amsmath}
\usepackage{amssymb}
\usepackage[pdffitwindow=false,bookmarks=true,pdfauthor={Matera et al.},colorlinks=true,citecolor=blue,linkcolor=blue]{hyperref}

\usepackage{natbib}
\usepackage[version=3]{mhchem}
\usepackage{mathtools}
\usepackage{mathrsfs}
\usepackage{amsxtra}
\usepackage{physics}
\usepackage{nccmath} 
\usepackage[]{algorithm2e, algorithmic}
\usepackage[version=3]{mhchem}

\newtagform {reaction} {\{} {\}}
\newcounter {save_equation}
\newcounter {reaction}

\newcolumntype{C}[1]{>{\centering\arraybackslash}p{#1}}

\begin{document}
\title{Multilevel Adaptive Sparse Grid Quadrature for Monte Carlo models}

\author{Sandra D\"opking}
\email{doepking@zedat.fu-berlin.de}
\affiliation{Institute for Mathematics, Freie Universit\"at Berlin, Arnimallee 6, D-14195 Berlin, Germany}

\author{Sebastian Matera}
\email{matera@math.fu-berlin.de}
\affiliation{Institute for Mathematics, Freie Universit\"at Berlin, Arnimallee 6, D-14195 Berlin, Germany}

\begin{abstract}
Many problems require to approximate an expected value by some kind of Monte Carlo (MC) sampling, e.g. molecular dynamics (MD) or simulation of stochastic reaction models (also termed kinetic Monte Carlo (kMC)). Often, we are furthermore interested in some integral of the MC model's output over the input parameters.  We present a Multilevel Adaptive Sparse Grid strategy for the numerical integration of such problems where the integrand is implicitly defined by a Monte Carlo model. In this approach, we exploit different levels of sampling accuracy in the Monte Carlo model to reduce the overall computational costs compared to a single level approach. Unlike existing approaches for Multilevel Numerical Quadrature, our approach is not based on a telescoping sum, but we rather utilize the intrinsic multilevel structure of the sparse grids and the employed locally supported, piecewise linear basis functions. Besides illustrative toy models, we demonstrate the methodology on a realistic kMC model for CO oxidation. We find significant savings compared to the single level approach - often orders of magnitude.  

\end{abstract}
\maketitle

\section{Introduction}
The need for the integration of high-dimensional functions arises in many scientific, engineering or socio-economic applications, e.g. in uncertainty quantification\cite{sullivan2015introduction}, finance\cite{griebel2010dimension}
or molecular simulations\cite{frenkel2001understanding}.
In most cases, the arising integrals can not analytically be solved and therefore need to be approximated by numerical quadrature. In many problem settings, the integrands are only implicitly given by a simulation algorithm and the evaluation of the integrand is the major source of computational costs. Additionally, simulations come usually with a finite (but tunable) error. In this study, we want to concentrate on integrands, which are defined by a Monte Carlo simulation and where the deviations of function value estimates results from the sampling error.
The need to parametrically integrate over the output of such Monte Carlo models arises, for instance, in the uncertainty quantification of stochastic reactivity models\cite{Alexanderian2014,DistributionGSA}. 
Another example is molecular dynamics simulations (MD), where the problem has multiple metastable sets, from which the MD can not escape, and an explicit quadrature over some degrees of freedom is used to overcome this problem\cite{chipot2007free}.\\

For numerical quadrature in high dimensions, a key prerequisite is that the numerical method does not suffer from the \textit{curse of dimensionality} \cite{bungartz2004sparse}; i.e. the exponential  increase of the computational costs with the dimensionality for a prescribed accuracy $\epsilon$.  Monte Carlo integration methods have a dimension independent but generally slow convergence, i.e. the error decreases as ${O}( N^{-1/2})$ when $N$ is the number function evaluations\cite{LEcuyer2005}. For medium high dimensionality, deterministic schemes like Quasi Monte Carlo (QMC) and Sparse Grids (SG) can instead achieve higher convergence rates, which only mildly depend on the dimensionality\cite{LEcuyer2005,gerstner1998numerical}.
For SG, a convergence of $O( N^{-r} (log N ) ^{r ( d-1)})$ can be achieved, where $r$ depends on the integrand and the employed basis functions. The foundations of SG were developed in the early nineties\cite{Zenger91,griebel1990parallelizable,bungartz1992dunne}, 
albeit similar ideas date back to the sixties\cite{Smo63}. As classical full tensor grids, Sparse Grids are based on products of one-dimensional interpolation and/or integration rules, but they omit certain \textit{higher order cross terms}, which do not contribute much to the accuracy. This largely removes the \textit{curse of dimensionality}, by which full tensor grids are affected. In this article, we concentrate on SG based on piecewise linear, locally supported basis functions, but the SG approach can be formulated using other basis functions such as polynomials\cite{nobile2008sparse}, wavelets\cite{webster2012adaptive} or prewavelets\cite{bungartz2004sparse}.

A key feature of SG is the possibility of adaptive refinement, i.e. choosing grid  points according to the function to be integrated and, thereby, further reducing
the amount of necessary points for numerical quadrature. The basis for adaptivity was already  established in the earliest papers on SG\cite{Zenger91,griebel1990parallelizable, bungartz1992dunne}. Adaptivity for high-dimensional interpolation and integration was first developed as the so-called generalized sparse grid (GSG) method \cite{griebel1998adaptive, Gerstner2003Griebel}, which, loosely speaking, adjusts the resolution for the different axes and has initiated a number of such dimension adaptive SG approaches\cite{GriebelHoltz10,Gerstner2003Griebel,bungartz2003multivariate}.
Later, local adaptivity for higher dimensional problems
 was introduced,\cite{ma2010adaptive,pflugerbook,webster2012adaptive},
which refines the grid close to rapid variations of the integrand. Different versions of combined local and dimension adaptivity have been proposed in recent years\cite{jakeman2011local,jakeman2012local,ma2010adaptive,stoyanov2013hierarchy,GSACo3O4}.

SG and especially adaptive Sparse Grids (ASG) can lower the number of possible expensive function evaluations. But, they do not provide a measure how accurately we have to estimate the function values by our simulations. The question arises which computational effort we have to spend and whether we should spend the same for every grid point.  Multilevel approaches tackle this problem by balancing the simulation with the quadrature error\cite{Giles2013}. The main idea behind multilevel methods is to utilize a hierarchy of numerical approximations of the underlying model. It turns out that often the most simulations can be performed with a low accuracy and the number of needed simulations decreases with increasing accuracy and computational cost. The most widespread family of multilevel methods for quadrature is Multilevel Monte Carlo (MLMC) with first appearance in the late nineties/early 2000s\cite{Heinrich98,Heinrich99,Heinrich2000}. MLMC can be applied to various problems, see ref. \cite{Giles2013} for an overview. Particularly, MLMC and its extension with QMC has attracted interest in the field of uncertainty quantification, e.g. for the treatment of partial differential equations with random input\cite{TeckentrupSchleichl, CliffeSchleichl}. In this context, the multilevel idea was also adapted for stochastic collocation with SG (MLSG) methods \cite{TeckentrupSC,Wyk2014}. 

In this study, we present a Multilevel Adaptive Sparse Grid (MLASG) approach for the parametric quadrature of Monte Carlo models, i.e. where the simulation error results from random sampling. In this case, the error of the model is not deterministic but stochastic. The telescoping sum approach of existing multilevel quadrature will then be of no benefit as the stochastic noise will simply produce random refinements. Therefore, we exploit the intrinsic multilevel structure of the Sparse Grid to guide the sampling of the Monte Carlo model. Using locally supported piecewise linear basis functions and local adaptivity, our findings show that the sampling effort can be halved with every refinement step without affecting the accuracy of the numerical discretization, which can save orders of magnitude in computational time.  

The manuscript is organized as follows: We outline the problem formulation in section \ref{sec:problem_formulation} where we introduce the general model and the relation between cost and error for function evaluations. We continue with a description of the employed SG method and the locally adaptive refinement strategy in section \ref{sec:sparse_grid}. In section \ref{sec:Multilevel}, we discuss the ideas, how this methodology can be exploited with multiple levels of sampling accuracy, and the final MLASG approach.
In section \ref{sec:results}, we present some illustrative examples to demonstrate the benefits of the MLASG approach. Besides some analytical test models, these examples also include a realistic stochastic model for catalytic CO oxidation.

\section{Parametric Monte Carlo models} 
\label{sec:problem_formulation}

For the rest of this manuscript, we will consider a function  $f: \Omega \rightarrow \mathbb{R}$ with $\Omega \in \mathbb{R}^D$, and we are interested in the integral
\begin{equation}\label{eq:Integral}
I = \int_{\Omega} f(x) dx^D.
\end{equation}
We restrict to those cases, where the problem can be formulated such that the domain $\Omega$ is a unit hypercube $\Omega = [0,1]^D$. We assume that the function $f$ is only implicitly
given by some simulation code, which takes $x$ as input parameters and  returns an approximation of $f(x)$. That is, we perform the integration over the input parameters and in a parametric Monte Carlo model the approximation of $f(x)$ is  done by some random sampling and the approximation errors are due to the finite number of samples. For a fixed $x$, our model produces samples from a probability distribution that is parametrically dependent on $x$, i.e. the output is a random variable $Y_x$. We assume that $f(x)$ is the expected value $\text{E} (Y_x)$ of $Y_x$ and is approximated by drawing $M_x$ samples and performing the statistical average, i.e.
\begin{equation}
f(x) = \text{E} (Y_x) \approx M_x^{-1} \sum_{i =1 }^{M_x} y_{x,i} =: \overline{Y_{xM}}
\end{equation}
where $y_{x,i}$ are the different samples according to the parameter $x$. Assuming that the simulation code provides independent samples and the variances and covariances of the estimates obey
\begin{equation}\label{eq:VarSample}
\text{Var} (\overline{Y_{xM}}) = C_x M_x^{-1}, ~\text{and} ~ \text{Cov} (\overline{Y_{xM}} ~\overline{Y_{x'M}} )=0 ~\text{for} ~ x\neq x'
\end{equation}
with $C_x = \text{Var}(Y_x)$. We assume that there exist an upper bound $C^*$ such that $C_x < C^*, \forall x \in \Omega$ . The variance of the approximation thus reduces with increasing sample size $M_x$, i.e. with increasing computational costs.  

In general, numerically approximating the integral \ref{eq:Integral}  with some kind of numerical quadrature, i.e.
\begin{equation}
I\approx I_N = \sum_{i=1}^N w_i f(x_i),
\end{equation} 
gives a sum of $w_i$-weighted  function evaluations $f(x_i)$ over $N$ samples . In the case of the Monte Carlo model the function evaluation have to be replaced by $\overline{Y_{x_i}}$ to solve the integral 
\begin{equation}\label{eq:ApproxQuadrature}
I_N \approx \overline{I_N} = \sum_{i=1}^N w_i \overline{Y_{x_iM}}.
\end{equation}
Since, the samples are drawn independently, the variance due to the noise in the function evaluations is given by 
\begin{equation}\label{eq:VarQuadrature} 
\text{Var} (\overline{I_N}) = \sum_{i=1}^N w_i^2 \text{Var} (\overline{Y_{x_i M}}) \leq C^* \sum_{i=1}^N w_i^2 
\end{equation}
and tends to zero for $N\rightarrow \infty$ since $\sum_{i=1}^N w_i^2 \rightarrow 0$ in that limit. In other words, we expect that $\overline{I_N}$ converges to the true integral for increasing numbers of grids points, irrespective of the values of the numbers of samples $M_{x_i}$ per grid point. This has consequences for the multilevel strategy, which we will discuss in section \ref{sec:Multilevel}.

\section{Sparse Grids }
\label{sec:sparse_grid}

 \begin{figure}
\begin{center} 
    \includegraphics[width=0.7\linewidth]{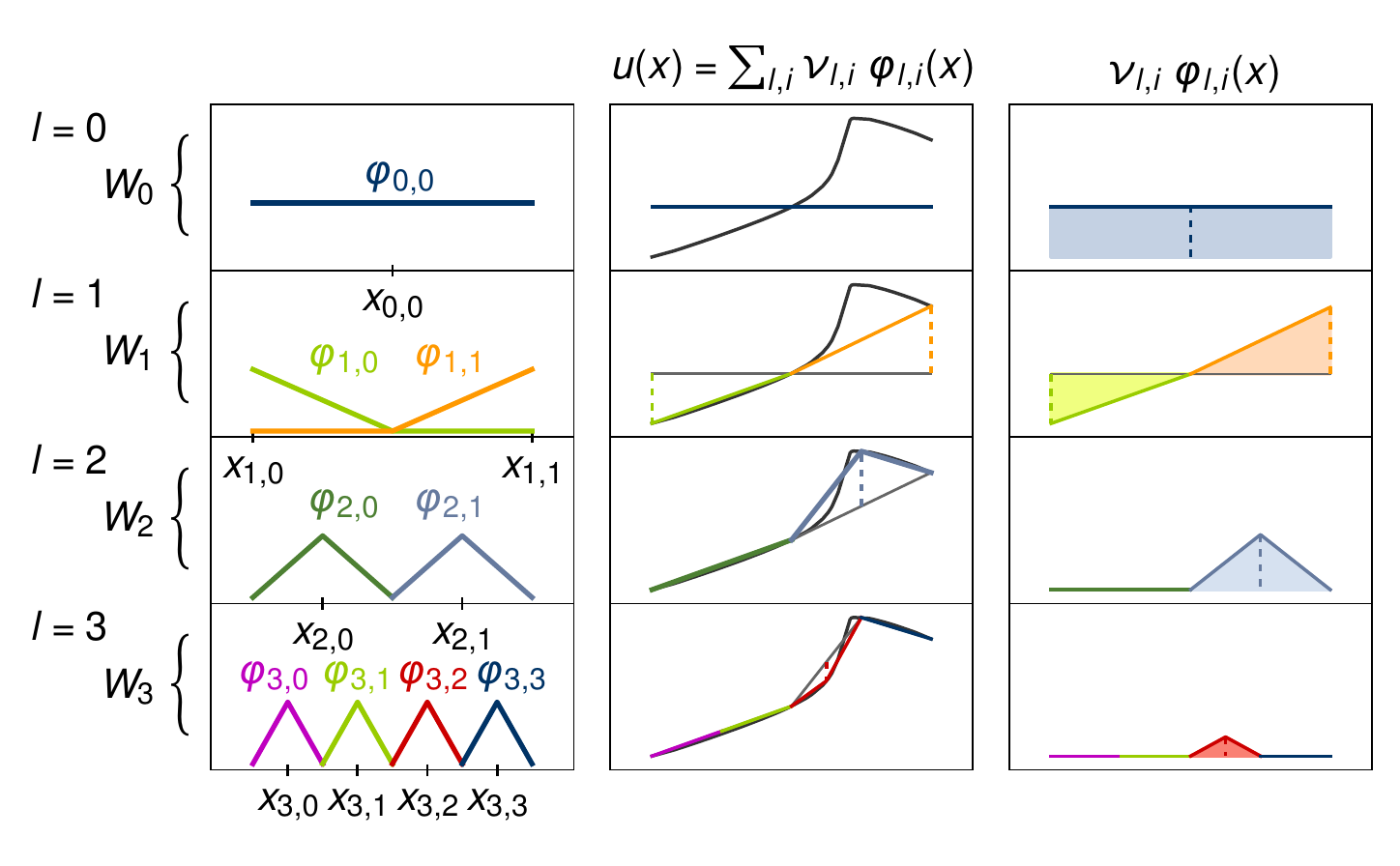}\\[-0.7cm]
\end{center}
\vspace{-0.5cm}
\caption{Illustration of the hierarchical interpolation in one dimension with
increasing level index. Left: Basis functions added at the corresponding levels $l$. Middle: The original function $f (x)$ and the hierarchical interpolation $\sum_{l\leq L ,i} \nu_{l,i} ~\varphi_{l,i}(x)$
Right: Illustration of the contribution of the basis functions added at each level $ \nu_{l,i}\varphi_{l,i}$. The contributions decrease with increasing level, and, at the finest level, only those basis functions close to the sharp non-linearity still have significant contributions to the interpolation.}
        \label{fig:hbf}
\end{figure}
 
In order to approximate the function $f$, the sparse grid method decomposes the space, in which the function $f$ lives, into contributions of hierarchical difference spaces\cite{bungartz2004sparse}.  These spaces are constructed from a product basis and depending on the characteristics of $f$, different kind of basis functions (BF) can be more or less efficient\cite{pflugerbook}. Here we concentrate on BFs which are products of hierarchical basis functions for piecewise linear approximation in one dimension\cite{klimke2005algorithm,ma2010adaptive}.

For a multivariate function $f: [0,1]^D \rightarrow \mathbb{R}$, the starting point is a space $V$ of univariate functions over one variable $x \in [0,1]$. We decompose $V$ into hierarchical difference spaces $W_l$, i.e. 
\begin{equation}
V =  W_0 \oplus W_1 \oplus ...\oplus W_{\infty} = \bigoplus_{l=0}^{\infty} W_l , 
\end{equation}
where $l$ is called the refinement level \cite{yserentant2006sparse,griebel1998adaptive}. Each of the finite dimensional spaces $W_l$ 
\begin{equation}
W_l =  \text{span} \left\{ \varphi_{l,i} | i = 0,...,m_l \right\}
\end{equation}
is spanned by $m_l$ basis functions $\varphi_{l,i}$. We choose the BFs\cite{klimke2005algorithm}
\begin{equation}\label{eq:BFunc1D}
\varphi_{l,i} (x) = 
\begin{cases} 
1 & \text{for} ~l = 0\\
\max(1-2^{l} |x-x_{l,i}|,0) & \text{for} ~l > 0,  \\ 
\end{cases}
\end{equation}
where $x_{0,0} = 0.5, x_{1,0/1} =0.5 \mp 0.5 ,x_{2,0/1} = 0.5 \mp 0.25$ and $x_{l,i} = x_{l-1, [i/2]}- (-1)^i/2^l$ with $i \in [0,2^{l-1}-1]$. Thus $m_l=2^l-1$ for $l=0/1$ and
$m_l=2^{l-1}-1$ else. 
 
These univariate BFs are depicted in the left panel of the figure \ref{fig:hbf}, from level $l=0$ to level $l=3$. The middle panel shows an example function and its piecewise linear interpolation $u(x)=\sum\limits_{l'\leq l, i} \nu_{l,i} \varphi_{l,i} (x) $ using all BFs up to this level $l$. This illustrates why $l$ is called the refinement level. For $l=0$, we only have the constant basis function, level $l=1$ adds those BFs which are needed for  a
piecewise linear approximation using the midpoint and the two ends of the interval. Every further level does the same; it
adds simply those BFs needed for a piecewise linear interpolation if we increase the resolution by placing additional nodes between the existing nodes \cite{pflugerbook, griebel1998adaptive}. The right panel of figure \ref{fig:hbf} shows the contributions $\nu_{l,i} \varphi_{l,i} (x)$ from the BFs of each level to the interpolant $u(x)$. If the function is well approximated between two nodes, additional BFs, with supports there, will have an only minor contribution and this can be exploited for an adaptive refinement as we will detail below. 

In the following, we will call those BFs from $l+1$ with the same support as a particular BF $ \varphi_{l,i}$ its
children and those from $l-1$ with the  same support its parent(s)\cite{jakeman2011local, jakeman2012local}. Correspondingly, all BFs from lower levels with the  same support are its ancestors and all BFs with the   same support from higher levels are its descendants. 

For a sparse grid approximation of a function $f:[0,1]^D \rightarrow \mathbb{R}$ in $D$ dimensions, we consider the spaces\cite{jakeman2012local}
\begin{equation}
V_{L}  = \bigoplus_{|\mathbf{l}|_1 \leq L} W_{\mathbf{l}}
\end{equation}  
where $\mathbf{l}=(l_0,...,l_D-1)$ is a multi-index denoting the level of refinement in each dimension $D$ and 
\begin{equation}
W_\mathbf{l} = \bigotimes \limits_{i=0}^{D-1} W_{l_i} = \text{span} \lbrace  \varphi_{\mathbf{l},\mathbf{i}} (\mathbf{x}):=\prod_{j= 0}^{D-1} \varphi_{l_j, i_j} (x_j) \, \text{with}\, i_j \in \{0,m_{l_j} \}\rbrace
\end{equation}
The difference to a full grid approximation is, that we use only those spaces $W_\mathbf{l}$ with $|\mathbf{l}|_1 < L$ instead of $|\mathbf{l}|_\infty < L$. Analogously to the one dimensional case, we will phrase two  product BFs $ \varphi_{\mathbf{l},\mathbf{i}}$ and $\varphi_{\mathbf{m},\mathbf{j}}$ parent and child if their absolute levels $|\mathbf{l}|_1$ and $|\mathbf{m}|_1$ differ by at most one and they have overlapping supports\cite{jakeman2011local, jakeman2012local}. If the absolute levels differ by more than one but the supports overlap, we will call them ancestor and descendant as in the one-dimensional case.

Using the Sparse Grid space $V_L$, we can approximates the function $f(\mathbf{x})$ by
\begin{equation}
f(\mathbf{x}) \approx u_{L}(\mathbf{x}) = \sum_{|\mathbf{l}|_1 \leq L ,\mathbf{i}} \nu_{\mathbf{l},\mathbf{i}} ~\varphi_{\mathbf{l},\mathbf{i}}(\mathbf{x})
\label{eq:Interpolation_SG}
\end{equation}
where the expansion coefficients $\nu_{\mathbf{l},\mathbf{i}}$ are also called hierarchical surplus. For numerical quadrature, we obtain the surplus by interpolating the $f$ at the nodes $\mathbf{x}_{\mathbf{l},\mathbf{i}}=(x_{l_0,i_0}, \ldots x_{l_{D-1},i_{D-1}})$, where $x_{l,i}$ are as for the definition of the univariate BF \ref{eq:BFunc1D}. The surpluses can be calculated according to the recursive formula\cite{bungartz2004sparse}
\begin{equation}\label{eq:Surplus}
\nu_{\mathbf{l},\mathbf{i}} = f(\mathbf{x}_{\mathbf{l}, \mathbf{i}}) - u_{|\mathbf{l}|_1 -1}(\mathbf{x}_{\mathbf{l},\mathbf{i}}). 
\end{equation}
The hierarchical surplus denote the hierarchical increments between two successive levels and corrects the coarser interpolation to the actual function value of $f$ at $\mathbf{x}_{\mathbf{l},\mathbf{i}}$. For sufficiently smooth functions, there is a bound for the hierarchical surplus, which tends to for $|\mathbf{l}|_1 \rightarrow \infty$\cite{bungartz2004sparse}. 

On the basis of the equation \ref{eq:Interpolation_SG} it is straightforward to define the quadrature. 
To approximate the integral of the function $f(\mathbf{x})$, we use the sparse grid interpolant \ref{eq:Interpolation_SG}
\begin{align}
\int_{\Omega} u_{L}(\mathbf{x}) dx^D &= \sum_{|\mathbf{l}|_1 \leq L,\mathbf{i}} \nu_{\mathbf{l},\mathbf{i}} ~ w_{\mathbf{l},\mathbf{i}}\\
\text{with} ~~ w_{\mathbf{l},\mathbf{i}}&= \int_{\Omega} \varphi_{\mathbf{l},\mathbf{i}} dx^D
\label{eq:quadrature}
\end{align}
with $w_{\mathbf{l},\mathbf{i}}$ as the weight of the basis function $\varphi_{\mathbf{l},\mathbf{i}}$. 

\subsection{Adaptivity}

While a non-adaptive sparse grids approach can tremendously reduce the number of necessary grid points, compared to a full grid with similar accuracy, it can still be that it adds grid points during refinement which do not improve much the accuracy of the interpolation or quadrature. This might for instance be because the function varies less in one or more directions than in the remaining ones or because strongly non-linear behavior appears only locally. Additional grid  points in subdomains, which are already well approximated, then add only minor improvements to the interpolation. Only adding grid points where needed would be desirable and the hierarchical children-parent relation enables such selective refinement.

In most cases, the important subdomains, where additional points are needed, are not a priori known, but have to be identified during the computational procedure. For this, locally adaptive refinement strategies have been developed, that attempt to reduce the number of grid points by refining only into those directions and in those subdomains, which can be associated with a high interpolation error. Practically, the exact local interpolation error is unknown, which implies an approximation of the local error. Adaptive strategies exploit here that, for BFs with a small contribution $\nu_{\mathbf{l},\mathbf{i}} ~\varphi_{\mathbf{l},\mathbf{i}}(\mathbf{x})$, usually also the descendants have a small contribution. For the adaptive refinement, some appropriate norm of the contribution $\nu_{\mathbf{l},\mathbf{i}} ~\varphi_{\mathbf{l},\mathbf{i}}(\mathbf{x})$ is employed as an indicator for the local error. 

Which norm to employ depends on the application. With our focus on numerical quadrature the 1-norm is the canonical choice and we therefore employ the indicator\cite{ma2010adaptive,jakeman2011local} 
\begin{equation}
\gamma_{\mathbf{l},\mathbf{i}} :=\left| \nu_{\mathbf{l},\mathbf{i}} \cdot w_{\mathbf{l},\mathbf{i}} \right|, 
\label{eq:error_criteria}
\end{equation}   
where $w_{\mathbf{l},\mathbf{i}} $ is the integration weight of $\varphi_{\mathbf{l},\mathbf{i}}$ which agrees with its 1-norm for the employed BFs. Often, the max-norm, i.e. the surplus $\nu_{\mathbf{l},\mathbf{i}}$,  is used as indicator\cite{pflugerbook,bungartz2004sparse}, 
but for locally rapid changes it decreases rather slowly and using a weighted surplus might still result in a sufficient accuracy, even when targeting at interpolation\cite{ma2010adaptive,jakeman2011local}. 

For the selective refinement we now calculate the indicator $\gamma_{\mathbf{l},\mathbf{i}}$ for all existing BFs with $|\mathbf{l}|_1=L$, where $L$ is the  highest level we currently have in the sparse grid. For all BFs, for which the indicator is above the predefined tolerance $tol$, we add the children and the  corresponding   grid points to the grid. So with every increase of the level from $L$ to $L+1$, those grid points and BFs which correspond to the set
\begin{equation}
\mathcal{A} = \left \{ (\mathbf{l}, \mathbf{i})| |\mathbf{l}|_1 = L+1 \wedge \max_{\mathbf{j}, \mathbf{m} \in p(\mathbf{l}, \mathbf{i})} (\left|\nu_{\mathbf{j},\mathbf{m}} \cdot w_{\mathbf{j},\mathbf{m}} \right|) > tol \right \},
\label{eq:set_of_points}
\end{equation}
are included, where $p(\mathbf{l},\mathbf{i})$ denotes the set of parents of the BF $\varphi_{\mathbf{l},\mathbf{i}}$. Additionally to the classical local refinement\cite{griebel1998adaptive, ma2010adaptive}, we include also all ancestors of the BFs defined by $\mathcal{A}$. Only for these two set of points, we evaluate the function and calculate the surpluses. Once we have  updated 
the sparse grid with this information, we repeat the procedure - now with a by one increased total level. We stop refining when the procedure finds no points for which $\gamma_{\mathbf{l},\mathbf{i}} \geq tol$.

By equation \ref{eq:Surplus}, this refinement strategy ensures that the surpluses of our adaptive strategy always agree with those from a full sparse grid, except for points which are not in the adaptive grid. We will exploit this property in the multilevel extension of the Adaptive Sparse Grids (ASG) approach.

\section{Multilevel adaptive quadrature}
\label{sec:Multilevel}
Besides the reduction of grid points by adaptive Sparse Grids (ASG), the adjustment of the accuracy of the approximate function evaluations is an additional way to improve the computational efficiency of the numerical quadrature. Reducing the overall computational effort for function evaluations is the goal of multilevel quadrature (MLQ) approaches, with Multilevel Monte Carlo (MLMC)\cite{giles2008multilevel} as the most prominent example, but also adaptive and non-adaptive Sparse Grids have been employed in this context\cite{TeckentrupSC,Wyk2014}. 

The idea is to employ a hierarchy of approximation levels $r \in [0,R]\subset \mathbb{N}$, where the accuracy, but also the computational costs, increase with $r$ and $R$ is the most accurate approximation which is still feasible. The different levels $r$ correspond to particular values for a numerical parameter, e.g.  a step size or a mesh spacing, if the function values of $f$ results from the solution of a (partial) differential equation. For MLQ, the integral of the most accurate approximation $I_R$ is rewritten using a telescoping sum
\begin{equation}
I_R = I_0 + \sum_{r>0}^R I_{r} - I_{r-1}= I_0 \sum_{r>0}^R \Delta I_{r}, ~ \text{with}~ I_r=\int_{\Omega} f_r(x) dx^D,
\end{equation}
where $f_r$ is the approximation to $f$ in level $r$. The integrals $I_0$ and $\Delta I_r$ are then numerically solved by independent quadrature rules. 
Under certain preliminaries\cite{CliffeSchleichl,Giles2013}, it turns out that the expensive integrals for high $r$ can be approximated with significantly less nodes than for small $r$ and 
the computational effort is significantly decreased compared to the direct numerical quadrature of $I_R$.

In our case, the accuracy is determined by sampling effort $M$ spend for an approximate function evaluation and a straightforward choice would be that $M$ should grow exponentially with $r$. However, since the errors for different grid points should be statistically independent, the discussion in section \ref{sec:problem_formulation} shows that $\Delta I_r$ should be zero.
Thus MLQ with non-adaptive quadratures reduces to a single numerical quadrature where the function values at the nodes can be approximated by a single sample from the Monte Carlo model.
A straightforward application of ASG to this problem would not work very well - if at all - because the large random noise for  a single sample prohibits a reasonably accurate estimation of the surplus. In consequence, we would simply get random refinements and the adaptivity would not pay off. 

We therefore do not follow the above outlined MLQ. Rather we will exploit the intrinsic ML structure of the sparse grids, where higher level BFs lead to a finer resolution. 
On the other hand, the contribution $\nu_{\mathbf{l},\mathbf{i}} w_{\mathbf{l},\mathbf{i}} $ of such BF also has potentially less impact on the value of the integral, because the weights decay exponentially with the level and surpluses should asymptotically decay, at least for functions with bounded mixed derivatives\cite{bungartz2004sparse}. We would expect that we can spend less effort for estimation of those function values.

The idea for a Multilevel Adaptive Sparse Grid (MLASG) approach is therefore to control the sampling effort for the estimates $\overline{Y_{\mathbf{l},\mathbf{i}}}$ of $f(x_{\mathbf{l},\mathbf{i}})$ such that the refinement strategy is not (much) affected. We achieve this by controlling
the variance of the (estimated) error indicator $\overline{\gamma}_{\mathbf{l},\mathbf{i}}$ \ref{eq:error_criteria}, such that
\begin{equation}\label{eq:MLASGbound}
\begin{split}
 \text{Var} (\overline{\gamma}_{\mathbf{l},\mathbf{i}}) \leq c \, tol^2\\
 \Leftrightarrow \text{Var}(\overline{\nu}_{\mathbf{l},\mathbf{i}}) \leq c\, tol^2 w_{\mathbf{l},\mathbf{i}}^{-2},
\end{split}
\end{equation}
with a user defined constant $c$. Since $w_{\mathbf{l},\mathbf{i}} \leq 2^{|\mathbf{l}|_1}$, we can increase the variance of the surplus by a factor four in each refinement step. Unfortunately, this does not directly translate to the estimates $\overline{Y_{\mathbf{l},\mathbf{i}}}$, whose accuracy we want to control. However,  by eq. \ref{eq:Surplus}, the surpluses are linear combinations of the estimates $\overline{Y_{\mathbf{l},\mathbf{i}}}$. As the later are statistically independent, we can then write eq. \ref{eq:MLASGbound} as 
\begin{equation}\label{eq:MLASGbound2}
 c\, w_{\mathbf{l},\mathbf{i}}^{-2} \, tol^2 \geq \sum\limits_{\mathbf{m},\mathbf{j}} A^2_{\mathbf{l},\mathbf{i};\mathbf{m},\mathbf{j}}  \cdot \text{Var}( \overline{{Y}_{\mathbf{m},\mathbf{j}}})
\end{equation}
where the matrix $A$ results from eq. \ref{eq:Surplus} and maps the function values to the surplus.

If we now add new points in a refinement step, a general procedure would work as follows: Let $M_{\mathbf{l},\mathbf{i}}$ be the sampling effort spend for the estimates $\overline{Y_{\mathbf{l},\mathbf{i}}}$ and $\Delta  M_{\mathbf{l},\mathbf{i}}$ the extra amount of sampling, which needs to be determined such that the inequality \ref{eq:MLASGbound2} is fulfilled. We would then minimize $\sum_{\mathbf{l},\mathbf{i}} \Delta  M_{\mathbf{l},\mathbf{i}}$ (the total sampling effort in this refinement step) subject to the constraint \ref{eq:MLASGbound2} utilizing eq. \ref{eq:VarSample}, or the corresponding inequality using the upper bound $C^*$. We would spend another $\Delta  M_{\mathbf{l},\mathbf{i}}$ for estimating $\overline{Y_{\mathbf{l},\mathbf{i}}}$, add the obtained estimates to the sparse grid and then proceed with the next refinement step.

Instead of solving the involved integer programming problem, we follow a different route and assume that we can find a constant $B$ such that choosing
\begin{equation}\label{eq:LevelVarRatio}
 \text{Var}( \overline{{Y}_{\mathbf{l},\mathbf{i}}}) = c \, tol^2 B^{|\mathbf{l}|_1}
\end{equation}
fulfills the ineq. \ref{eq:MLASGbound2}. Under this assumption, the inequality \ref{eq:MLASGbound2} is fulfilled
 in every refinement step if 
\begin{equation}\label{eq:MLASGbound3}
 \sum\limits_{\mathbf{m},\mathbf{j}} A^2_{\mathbf{l},\mathbf{i};\mathbf{m},\mathbf{j}} B^{-|\mathbf{m}|_1}  \leq 4^{|\mathbf{l}|_1} 
\end{equation}
where we used that $w_{\mathbf{l},\mathbf{i}} \leq 2^{|\mathbf{l}|_1}$. We can further exploit that $A^2_{\mathbf{l},\mathbf{i};\mathbf{m},\mathbf{j}}=0$ for $|\mathbf{m}|_1 > |\mathbf{l}|_1$ by equation \ref{eq:Surplus} and that for fixed $(\mathbf{l},\mathbf{i})$ we have the same coefficients $A^2_{\mathbf{l},\mathbf{i};\mathbf{m},\mathbf{j}}$ in our adaptive strategy as if we would employ no adaptive refinement. We can now test different values of $B$ for different dimensionalities and up to a certain level using the $A^2_{\mathbf{l},\mathbf{i};\mathbf{m},\mathbf{j}}$ from non-adaptive SG. For the range of dimensions and maximum levels in this study, we found that the inequality \ref{eq:MLASGbound3} is fulfilled
for choosing $B=2$. 

What 
is left is a good choice for the constant $c$. In our experiments, we found that $c=1$ is fully sufficient to achieve a well working adaptive refinement. Values below do not improve much the final accuracy and significantly larger values lead to the afore mentioned random refinements. We explain this rather large value by our choice for the level to variance ratio \ref{eq:LevelVarRatio}, which is likely rather conservative. Therefore, the surpluses might actually have a significant lower variance than $tol^2$. Further, as long as the surplus is much larger than $tol$ a corresponding variance is sufficient to make false refinement or non-refinement very unlikely. If the surplus is in the order of $tol$, false refinement or non-refinement become likely. However, the surplus is only a rough error indicator. Possible errors due to this indicator have likely a similar impact as if we erroneously refine or not some BFs due to the noise of order $tol$.

The discussion on the general MLQ and the consequences, when it is applied to parametric Monte Carlo models, indicate that MLASG is not always superior for quadrature compared to the simple non-adaptive case. In the non-adaptive case, we can simply run all simulations at the lowest accuracy and would rather spend many grid points. For MLASG, the first point must be run at rather high accuracy and high computational costs. So the benefit of MLASG for quadrature only appears, when the adaptive refinements can sort out large numbers of grid points. In other words, we need to achieve a higher refinement level with MLASG than with the brute-force approach for the same computational costs. On the other hand, the error on the interpolation due to random noise in the samples seems to be always below the actual interpolation error of the ASG for a choice of $c\leq 1$. So MLASG seems to be able to also accelerate the construction of surrogates for Monte Carlo models, i.e. a similar setting as it has originally  been treated by Heinrich in his pioneering work on MLMC\cite{Heinrich98,Heinrich99,10.1007/3-540-45346-6_5}.
This property shall be investigated in future studies as such surrogates have a wide range of applications.

\section{Results}\label{sec:results}

\subsection{Test-model}
\label{subsec:TestModel}
As a first test case, we consider the hypercube of $[-0.5,0.5]^D$  and the function 
\begin{equation}
f ({x})= g\left(\sqrt{\sum_i^d (x_i +0.5)^2}\right) , 
\end{equation}
with 
\begin{equation}
g(r) =
\begin{cases}
g_1(r)=10  (\exp(\dfrac{-r+0.35}{0.086}) +1)^{-1} ~, \text{for} ~ r < 0.6\\ 
g_2(r)=0.005^r  g_1(0.6) \dfrac{1}{(0.005^{0.6})} ~, \text{for} ~ r \geq 0.6
\label{eq:toymodel}
\end{cases}.
\end{equation}
By construction the function $f$ has a kink for $\sqrt{\sum_i^d (x_i +0.5)^2}=0.6$ and it is shown in the left panel of figure \ref{fig:2dGrids} 
for $D=2$. 
We mimic the output of a Monte Carlo model by adding normally distributed, zero mean noise to the function evaluations, i.e. 
\begin{equation}
Y_{\mathbf{l},\mathbf{i}} = f(x_{\mathbf{l},\mathbf{i}}) + s_{\mathbf{l},\mathbf{i}},~ \text{with}~ \mathrm{Var}(s_{\mathbf{l},\mathbf{i}})=c\,tol^2\, 2^{|\mathbf{l}|}
\label{eq:noise_f}
\end{equation}
where $s_{\mathbf{l},\mathbf{i}}$ is the noise and $Y_{\mathbf{l},\mathbf{i}}$ are the values used for the sparse grid construction. In real MC models, the computational cost for obtaining $Y_{\mathbf{l},\mathbf{i}}$ is related to the variance by eq. \ref{eq:VarSample}. We therefore assign a cost ratio of $2^{-|\mathbf{l}|}$ to the sample $Y_{\mathbf{l},\mathbf{i}}$ and, in this way, mimic a real parametric Monte Carlo model.

\subsubsection{2D-case}
\begin{figure}
\begin{center} 
\includegraphics[width=0.49\linewidth]{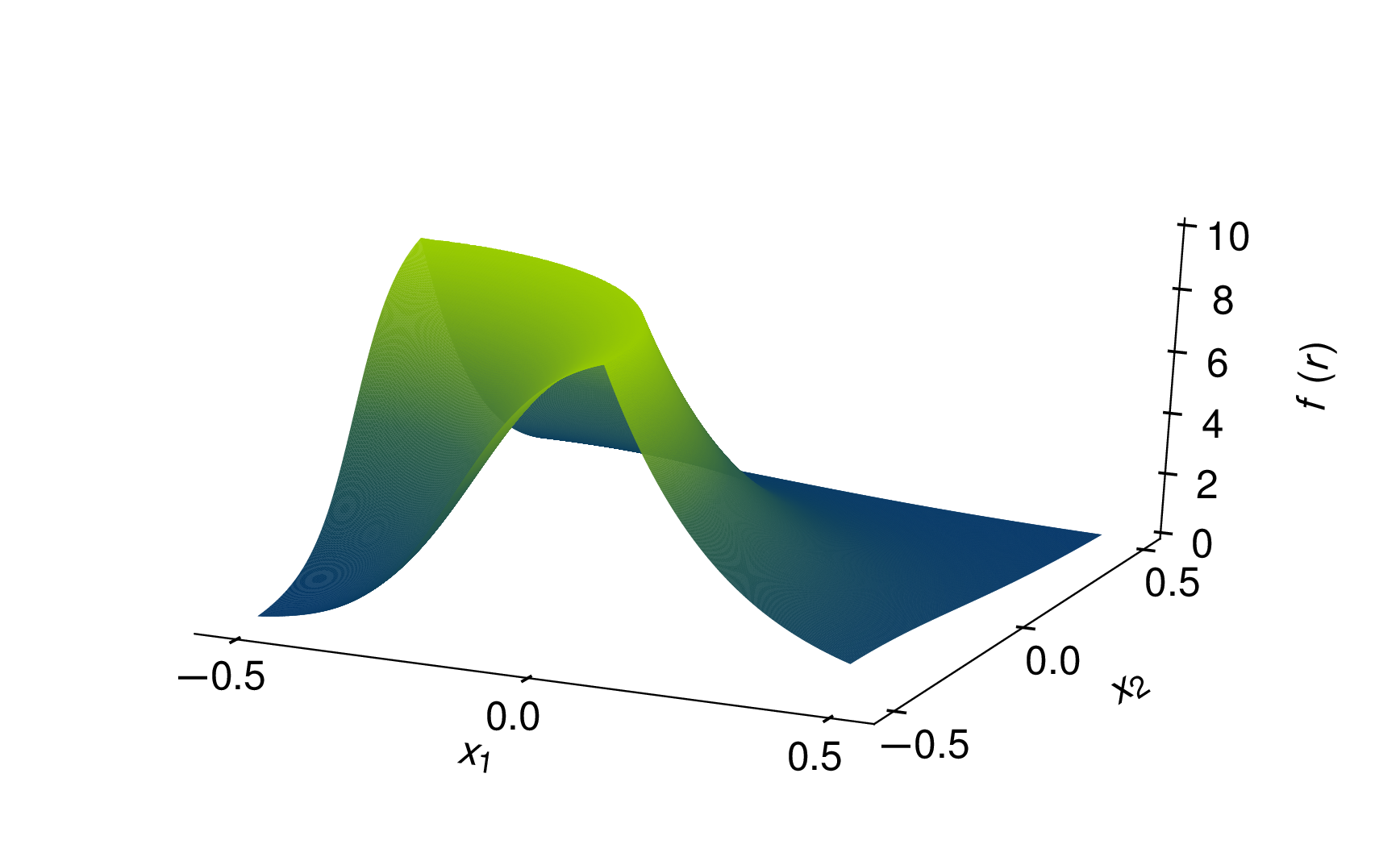}
    \includegraphics[width=0.49\linewidth]{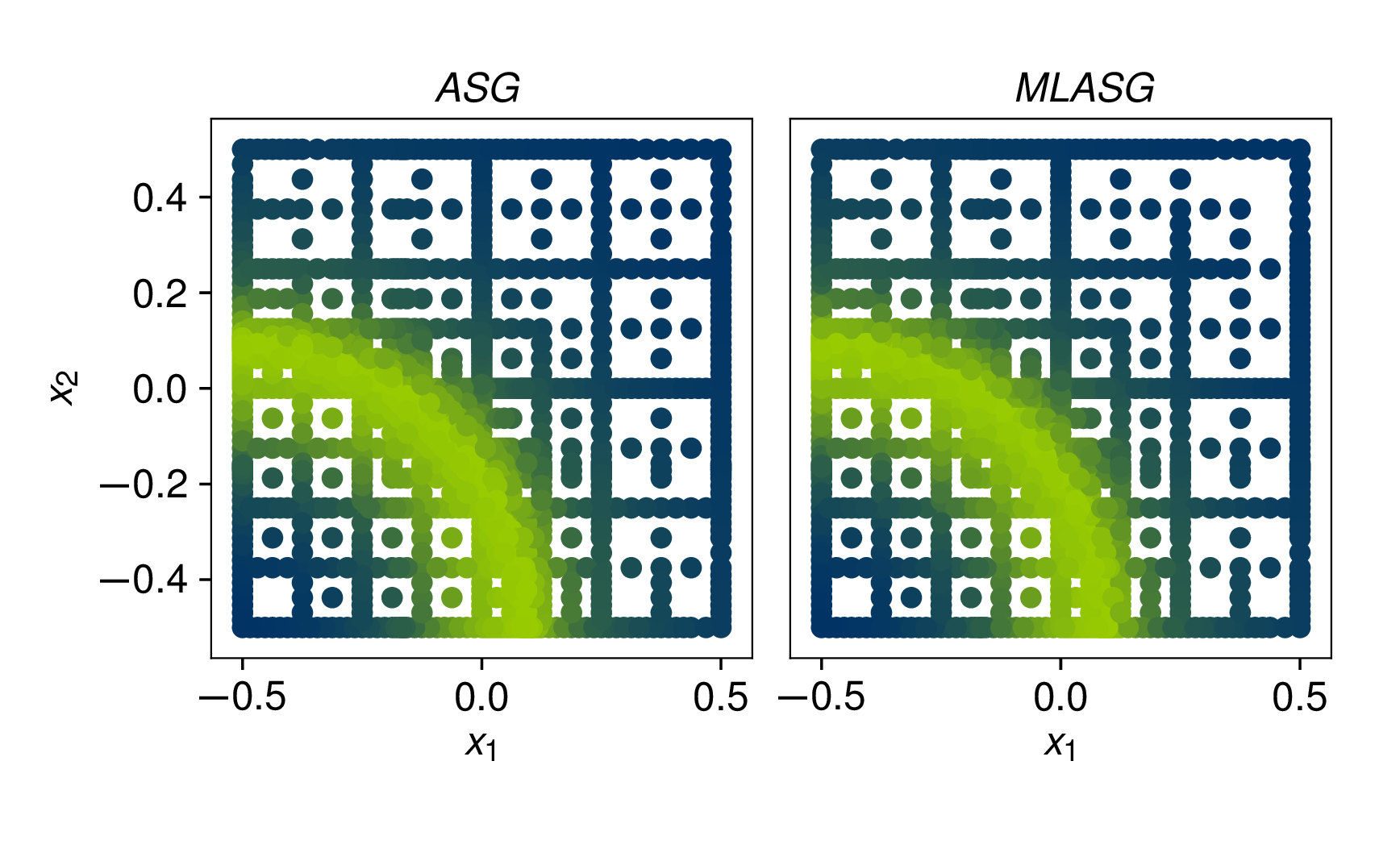}
\end{center}
\vspace{-0.5cm}
\caption{\textit{Left panel}: A representation of the test function \ref{eq:toymodel} in a 2-dimensional space with a kink defined on the radial coordinate $r=0.6$. \textit{Right panel}:Adaptive grids for ASG (left) and MLASG (right) for the test problem \ref{eq:toymodel} in the 2-dimensional case employing the tolerance $tol=10^{-4}$ and $c=1$ (compare eq. \ref{eq:MLASGbound}) }
        \label{fig:2dGrids}
\end{figure}

\begin{figure}
\begin{center} 
    \includegraphics[width=0.49\linewidth]{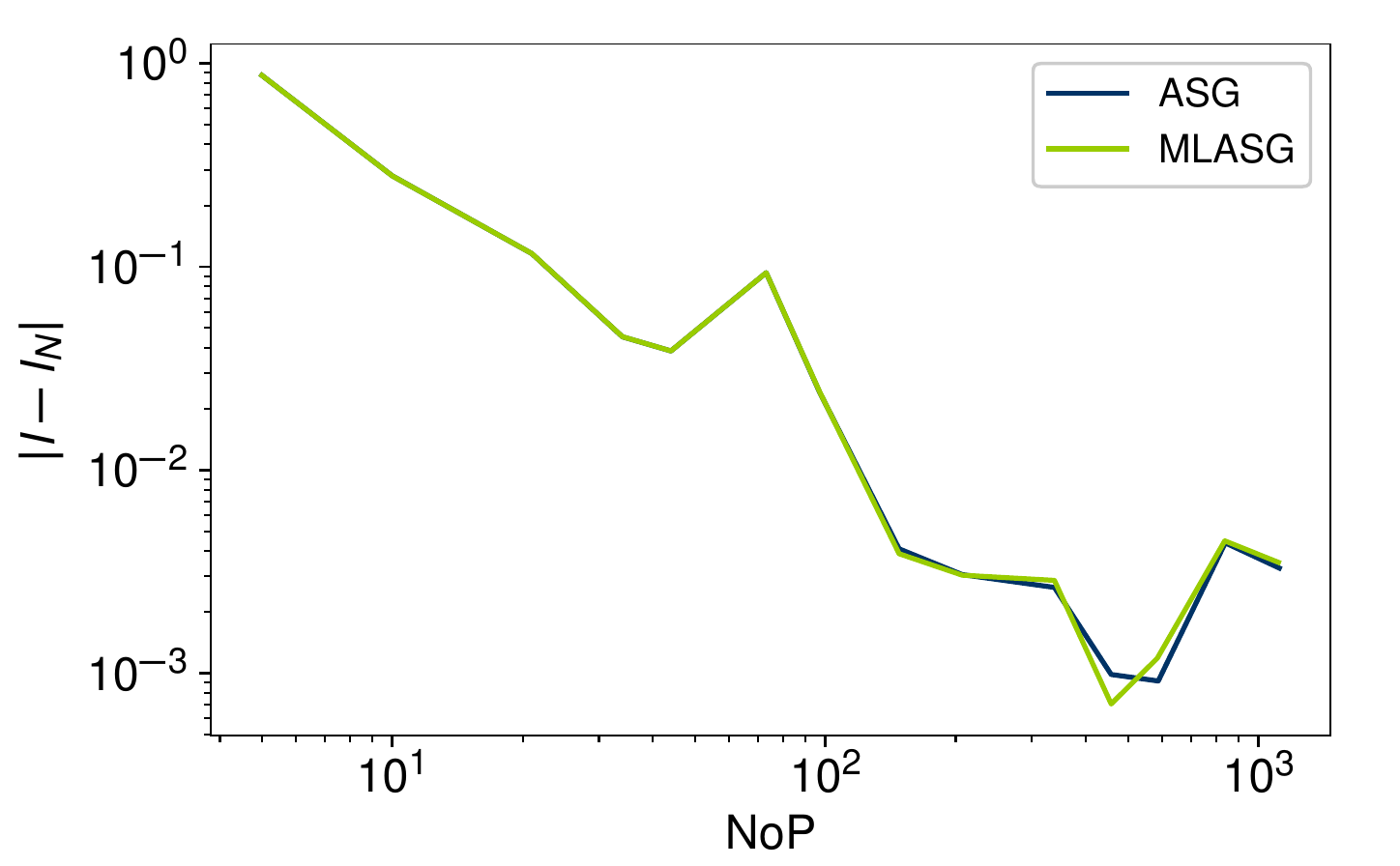}
   \includegraphics[width=0.49\linewidth]{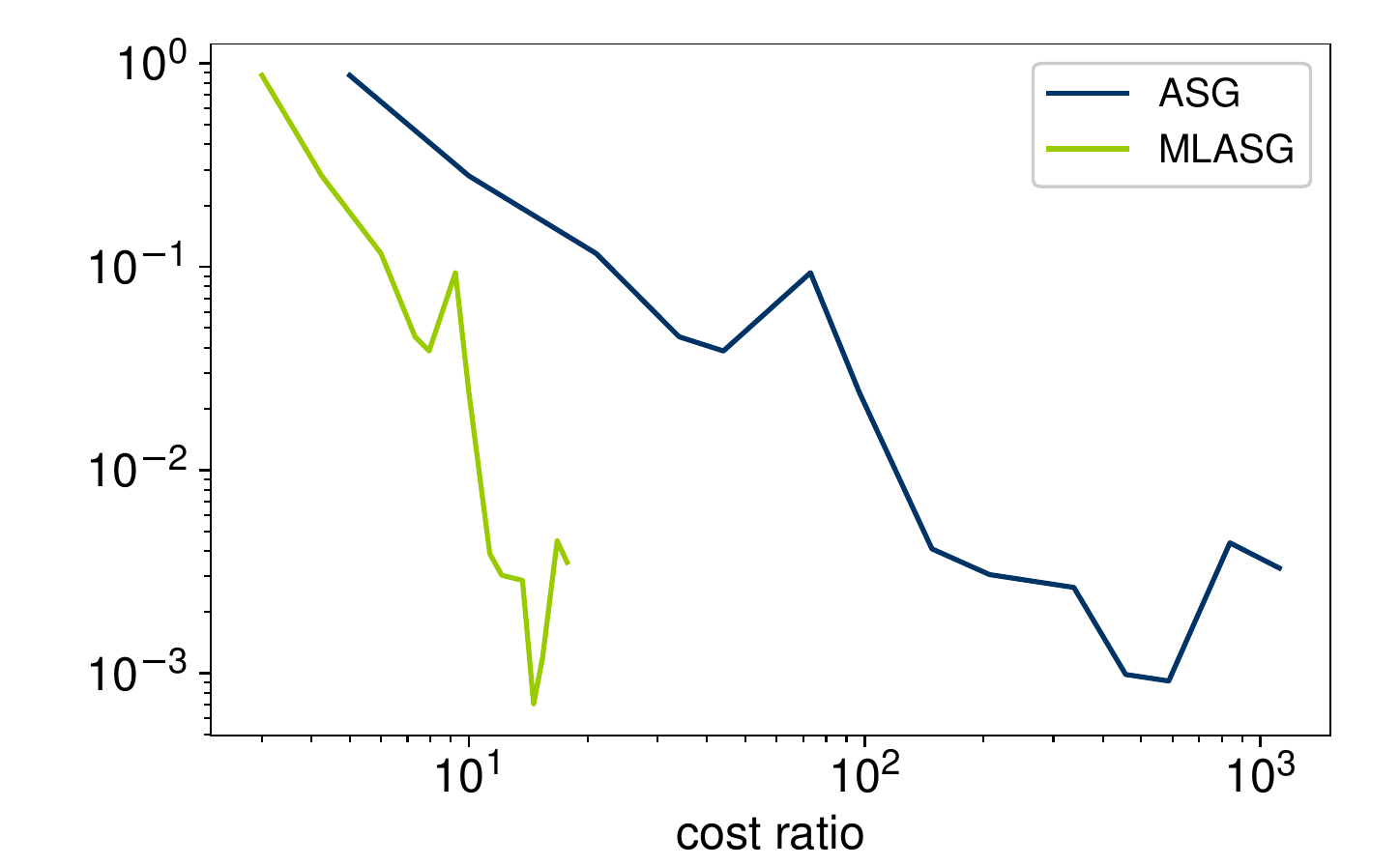}
\end{center}
\vspace{-0.5cm}
\caption{Comparison the quadrature error $|I-I_N|$ during the refinement process of ASG and MLASG for the test model \ref{eq:toymodel} in the 2-dimensional case. For this a tolerance $tol=10^{-4}$ and $c=1$ (compare eq. \ref{eq:MLASGbound}) are employed. 
\textit{Left panel}: Quadrature error as function of the number of grid points (NoP).
\textit{Right panel}: Quadrature error as function of the expected computational cost (in multiples of the cost for the first point).}%
        \label{fig:2DIntegrals}
\end{figure}

In order to visualize the effect of the MLASG on the adaptive refinement, we first consider a dimensionality of $D=2$ with a threshold
$tol=10^{-4}$. We employ $c=1$ for the constant in eq. \ref{eq:MLASGbound}, i.e. the standard deviation $\sigma_0$ for the first point (level $L=0$) equals the threshold $tol$. The adaptive grids for the standard single level Adaptive Sparse Grids (ASG) and the Multilevel Adaptive Sparse Grid (MLASG) approach are shown in the right panel of fig. \ref{fig:2dGrids} on the left and on the right, respectively. The MLASG produces almost the same adaptive grid 
as the standard procedure refining the grid close to kink of $f$. This good agreement of the MLASG with the reference ASG translates to the accuracy of the numerical quadrature. Both approaches show the essentially same convergence behavior during the refinement process. This can be seen from the left panel in figure \ref{fig:2DIntegrals}, where we plot the quadrature error during the refinement against the number of grid points (NoP). Both curves essentially lay on top of each other and deviations are only visible, when the accuracy is already very low. Even then the differences are such, that it is not possible to judge whether ASG or MLASG is better. This supports the discussion from the previous section \ref{sec:Multilevel} on why the choice $c=1$ is good enough.

The advantage of the MLASG over the standard ASG can be seen from the right panel in figure \ref{fig:2DIntegrals}. Again, we show the quadrature error during refinement, but now in dependence of the expected computational cost, which would be spend if our toy model would be a real Monte Carlo model. While initially, at the level $L=1$, the  evaluation
costs are very similar, the MLASG converges much more rapidly with respect to the costs. After termination, we have spend two orders of magnitude less resources in MLASG than in ASG.

\begin{figure}
\begin{center} 
    \includegraphics[width=0.49\linewidth]{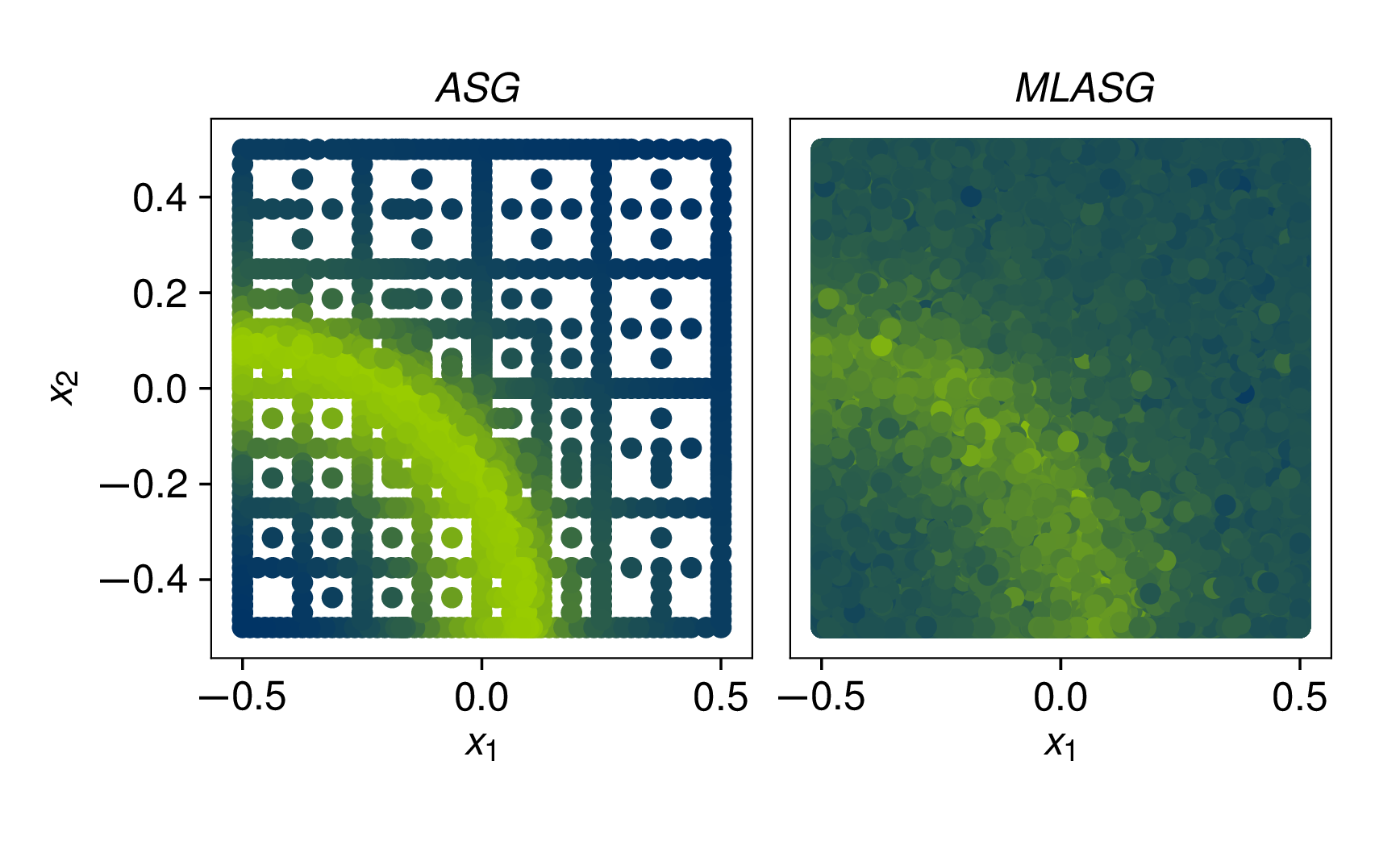}
    \includegraphics[width=0.49\linewidth]{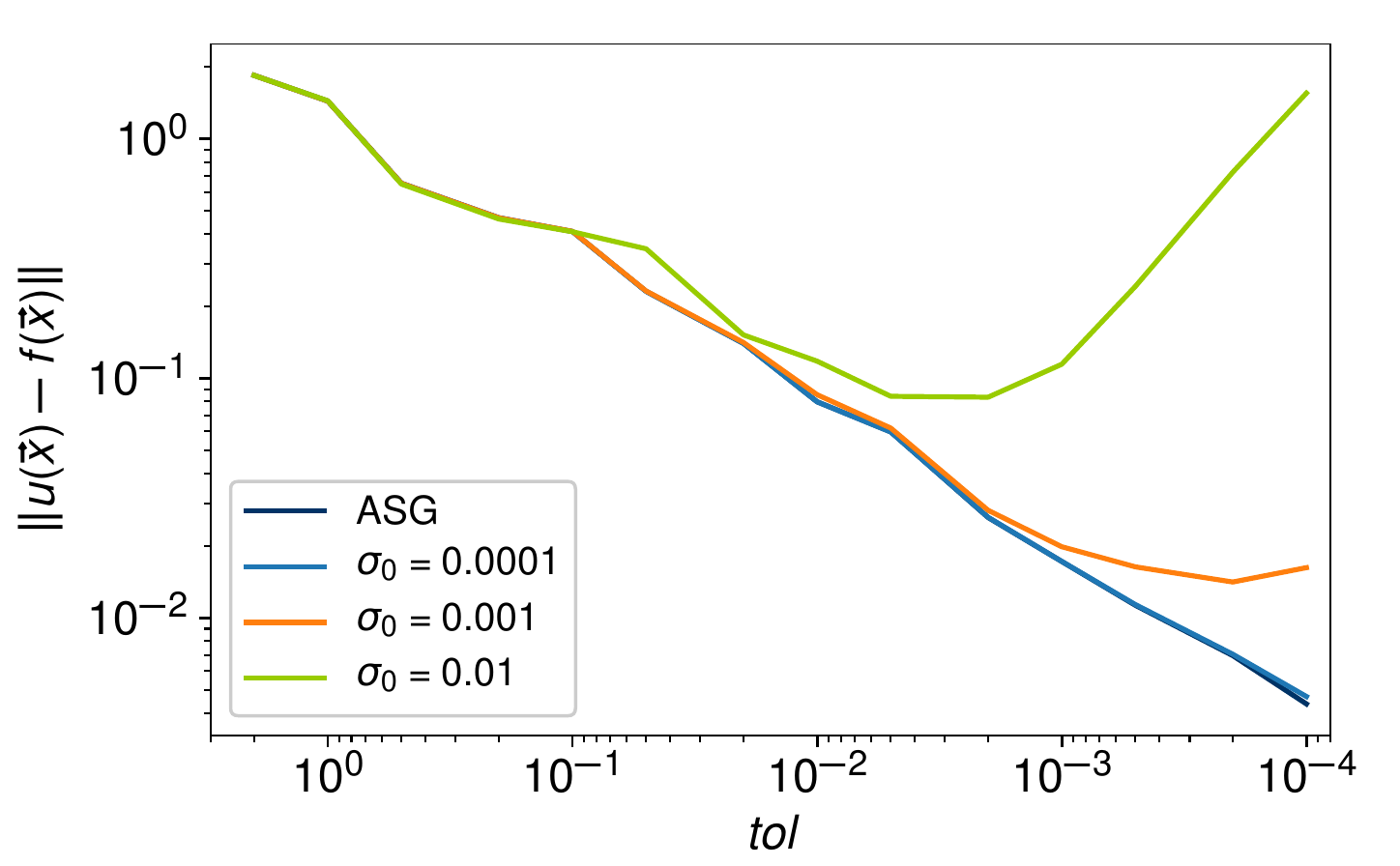}
\end{center}
\vspace{-0.5cm}
\caption{\textit{Left panel}: Comparison of the adaptive grids for $tol=10^{-4}$ for ASG without noise and MLASG with an initial standard deviation $\sigma_0=10^{-2}$, i.e. $c=100$.
\textit{Right panel}: 1-norm of the interpolation error for different initial noise values $\sigma_0 = [10^{-2}, 10^{-4}]$ compared to the the ASG method without noise.}
        \label{fig:AR_sigma_2D}
\end{figure}

We now turn to a more detailed investigation of the choice of the constant $c$ controlling the ratio between the sampling variance on the refinement tolerance $tol$. The left panel in figure \ref{fig:AR_sigma_2D}, shows the adaptive grids for the ASG using the exact function values and MLASG with an initial standard deviation $\sigma_0=10^{-2}$. For both, we employed $tol=10^{-4}$, i.e. the only difference to the previous example is that we now have chosen $c=100$ instead of $c=1$ for the MLASG. The effect is quite large. While for $c=1$ both grids essentially agreed (fig. \ref{fig:2dGrids}) , the MLASG grid does not follow the characteristics of the function and grid points have essentially been placed everywhere in the domain. The reason for that is obviously the fact that the noise in the refinement criterion is a factor hundred larger than the threshold. Therefore, false decision in the adaptive selection are quite likely also for points, which would have a comparatively large surplus. 

In order to analyze the effect, we have performed  MLASG for varying values of the threshold $tol \in [10^{-4}, 2]$ and for three values of the initial standard deviation $\sigma_0 \in\{ 10^{-2},10^{-3},10^{-4}\}$. The right panel in fig. \ref{fig:AR_sigma_2D} shows the  error of these three cases as a function of the threshold $tol$. For comparison  we added the ideal ASG case.
Albeit we focus on quadrature, we display the 1-norm of the difference between the true function $f$ and the respective interpolant $u(x)$, because the effect of the choice of $c$ is better visible in this measure. Starting at the largest values for $tol$, all three curves follow the ideal ASG curve until $tol$ becomes close to the respective value of $\sigma_0$, i.e. when $c\approx 1$. For values below, the curves start to flatten
and from a certain point on the error does actually increase again. Thus, choosing $c$ (much) smaller than one does not significantly improve the accuracy. For much larger values, we loose the benefits of adaptive refinement and, at least, the interpolation accuracy.  

\subsubsection{7D-case}
\begin{figure}
\begin{center} 
    \includegraphics[width=0.49\linewidth]{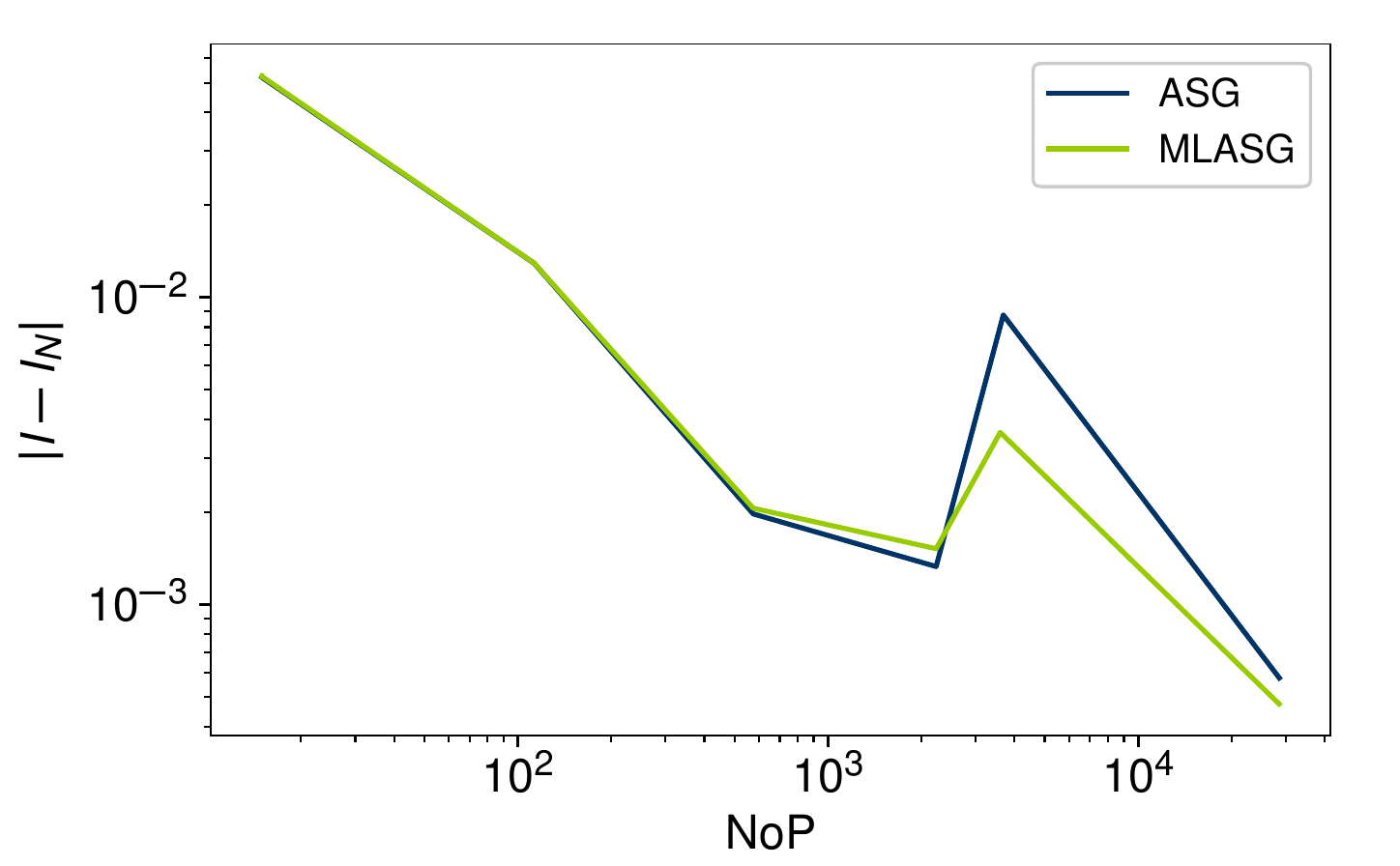}
   \includegraphics[width=0.49\linewidth]{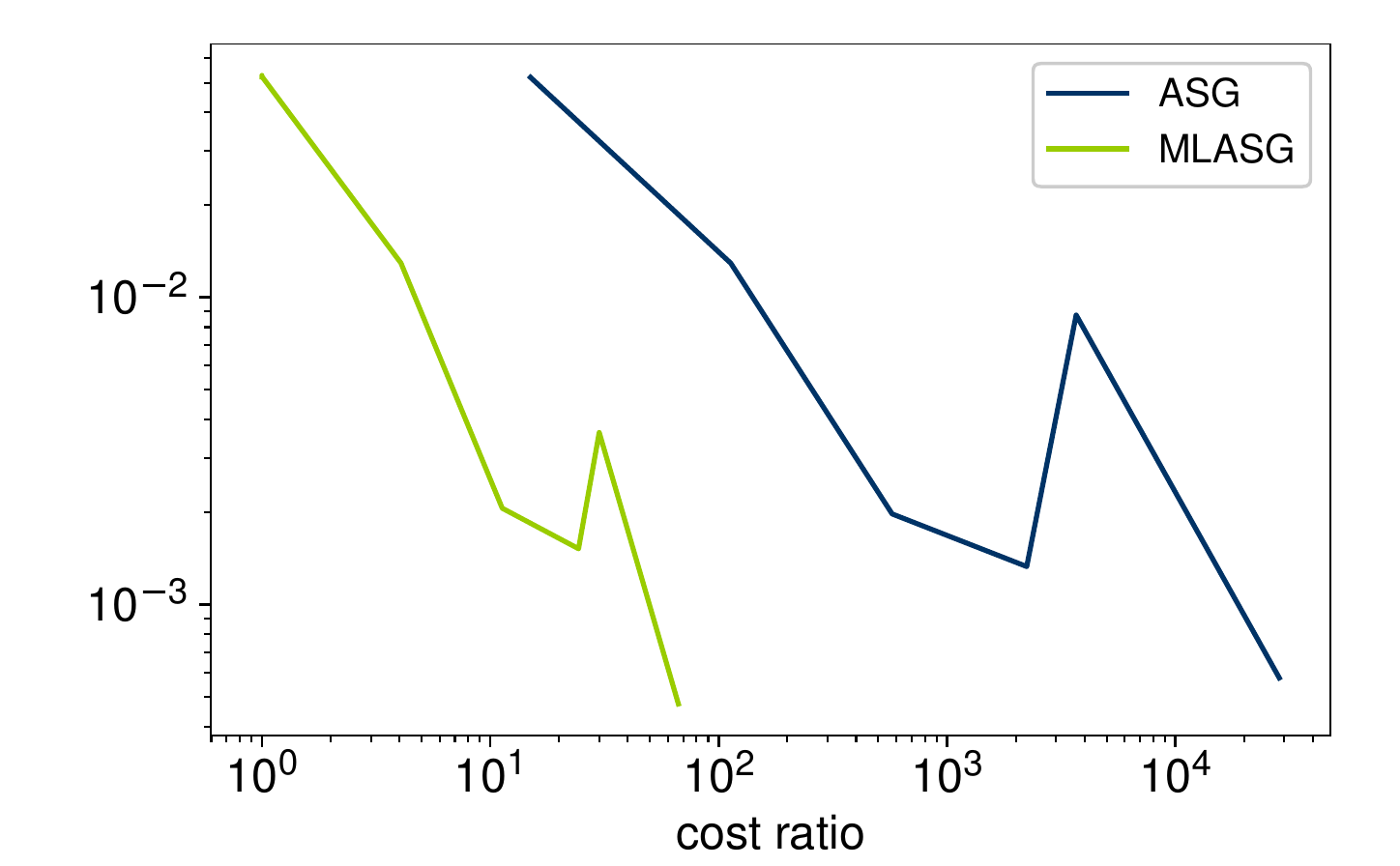}
\end{center}
\vspace{-0.5cm}
\caption{Comparison the quadrature error $|I-I_N|$ during the refinement process of ASG and MLASG for the test model \ref{eq:toymodel} in the 7-dimensional case. For this a tolerance $tol=10^{-4}$ and $c=1$ (compare eq. \ref{eq:MLASGbound}) are employed. 
\textit{Upper panel}: Quadrature error as function of the number of grid points.
\textit{Lower panel}: Quadrature error as function of the expected computational cost (in multiples of the cost for the first point).}%
        \label{fig:7DIntegrals}
\end{figure}
As a higher dimensional example, we consider the test function \ref{eq:toymodel}, but now the 7-dimensional case. Figure \ref{fig:7DIntegrals} shows the equivalent to fig. \ref{fig:2DIntegrals} for this case, i.e. the quadrature error during refinement as function of the number of grid points (left panel) and the cost ratio (right panel). As before, we employed a refinement tolerance of $tol=10^{-4}$ and $c=1$. Both, ASG and MLASG, show very similar behavior with the number of grid points and we would conclude that, as in the previous example, the adaptive refinement is not much affected by the increasing variance during refinement in MLASG. As in the 2D example, this increase with the refinement level has the effect that MLASG can save orders of magnitude in computational effort for the function evaluation as it is visible from the right panel in fig. \ref{fig:7DIntegrals}. In the 7-dimensional case, the savings are even larger, almost three orders of magnitude.

\subsection{CO oxidation model}
\label{subsec:kMCModel}
As a realistic example, we turn to a stochastic model for the CO oxidation on a RuO$_2$(110) heterogeneous catalyst\cite{1DGelss,lorenzi2017local} and address its parametric uncertainty. 
It is a reduced version of the original quantum chemistry based model \cite{ReuterSchefflerRuO2} and exploits that the chemistry is mainly controlled by the so-called cus adsorption sites \cite{meskine2009examination,hoffmann2016,DistributionGSA}. It describes the chemical kinetics of the CO oxidation as Markov jump process on a chain of these sites on the catalytic surface and every jump corresponds to one of the allowed elementary reactions. The elementary reactions are the adsorption and desorption of CO and oxygen,  diffusion from on site to a  neighboring site of adsorbed CO and oxygen and the formation gaseous CO$2$. 

In the model, each site, enumerated by its position $l$ in the chain, can be in one of three states: I) empty ($\ce{e_l}$), II) CO covered ($\ce{CO_l}$), or III) oxygen covered ($\ce{O_l}$). The state of the chain, i.e. the vector carrying the states of the individual sites, is changed by the elementary processes
\begin{align*}
 	\ce{e_l &->[k^{\rm ads}_{\rm CO}] CO_l},\\
	\ce{e_l + e_{l+1} &->[k^{\rm ads}_{\rm O_2}] O_l + O_{l+1}},\\
  	\ce{CO_l &->[k^{\rm des}_{\rm CO}] e_l},\\
	\ce{O_l + O_{l+1} &->[k^{\rm des}_{\rm O_2}] e_l + e_{l+1}},\\
  	\ce{CO_l + e_{l+1}&->[k^{\rm diff}_{\rm CO}] e_l+ CO_{l+1} },\\
  	\ce{e_l + CO_{l+1}&->[k^{\rm diff}_{\rm CO}] CO_l+ e_{l+1} },\\
  	\ce{O_l + e_{l+1}&->[k^{\rm diff}_{\rm O}] e_l+ O_{l+1} },\\
  	\ce{e_l + O_{l+1}&->[k^{\rm diff}_{\rm O}] O_l+ e_{l+1} },\\
  	\ce{CO_l + O_{l+1}&->[k^{\rm reac}] e_l+ e_{l+1} },\\
  	\ce{O_l + CO_{l+1}&->[k^{\rm reac}] e_l+ e_{l+1} },\\
\end{align*}
where $k^{\rm ads}_{\rm CO}$ etc. are the rates, with which the respective process occurs. For further details on the model, we refer to the references  \cite{1DGelss,lorenzi2017local}.

By the dependence on the rates, the model depends on seven independent parameters and we choose following parametrization: I) the equilibrium constant of adsorption of CO ($K_{\text{CO}}$) and oxygen ($K_{\text{O}_2}$), II) the rates  $k^{\rm ads}_{\rm CO}$ and $k^{\rm ads}_{\rm O_2}$ for these two reactions, III) the  rates $k^{\rm diff}_{\rm CO}$ and $k^{\rm diff}_{\rm O}$ for diffusion of adsorbed species and IV) the rate $k^{\rm reac}$ for the formation of CO$_2$. The rates for desorption of CO and O$_2$, $k^{\rm des}_{\rm CO}$ and $k^{\rm des}_{\rm O_2}$, are implicitly determined by the relation $ k^{\rm des}_{\rm CO/ O_2}=k^{\rm ads}_{\rm CO/ O_2}K^{-1}_{\rm CO/ O_2}$. 

The input parameters of the model have been obtained from Density Functional Theory (DFT). As the true density functional is only approximately known, the parameters carry some finite error. An uncertainty analysis requires to integrate over these parameters using a probability weight, which reflects the uncertainty in the parameters. For the present model case, we assume that the logarithms of the input parameters are independently and uniformly distributed in a certain range of values. This reflects that the uncertainties are dominated by the errors in binding energies and barriers, which enter the input exponentially\cite{ReuterSchefflerRuO2}. The employed ranges as well as the default values are provided in the table \ref{tab:Rate_constants}. 

\begin{table} 

\caption{\label{tab:Rate_constants} List of input parameters  \\ }

\begin{tabular}{l|c|c}

Name & Default value [1/s] & Range [1/s] \\[0.1cm]



\toprule \\[-0.2cm]

$K_{\text{CO}}$ & $\frac{2.0}{9.2} \times 10^{2}$ & $\frac{2.0}{9.2} \times 10^{0}$ - $\frac{2.0}{9.2} \times 10^{4}$  \\

$K_{\text{O}_2}$ &$\frac{9.7}{2.8} \times 10^{6}$ & $\frac{9.7}{2.8} \times 10^{4}$ - $\frac{9.7}{2.8} \times 10^{8}$ \\

$k^{\rm ads}_{\rm CO}$  & $2.0 \times 10^{8}$ & $1.0 \times 10^{8}$ - $4.0 \times 10^{8}$ \\

$k^{\rm ads}_{\rm O_2}$ & $9.7 \times 10^{7}$ & $4.85 \times 10^{7}$ - $1.94 \times 10^{8}$ \\

$k^{\rm diff}_{\rm CO}$ & $5.0 \times 10^{-1}$ &   $5.0\times 10^{-3} $ - $5.0 \times 10^{1} $\\

$k^{\rm diff}_{\rm O_2}$ & $6.6 \times 10^{-2} $ & $6.6 \times 10^{-4} $ - $6.6 \times 10^{0}$ \\

 $k^{\rm reac}$ & $1.7 \times 10^{5} $ &  $1.7 \times 10^{3} $ - $1.7 \times 10^{7} $\\

\end{tabular}

\label{Tab.: rate_constant_CO}

\end{table}


As an example for the uncertainty analysis, we consider the parameter averaged stationary CO coverage, where the coverage is the (process averaged) concentration of a certain type of adsorbate in units of number of molecules per site. The coverage is therefore always between $0$ and $1$. The required integral involving the outlined parameter distribution and  the parameter dependent CO coverage (as the output of the kMC simulation) can easily be brought into the form \ref{eq:Integral} by a coordinate transformation. The employed parameter ranges are now such that the coverage is close to one of these two values in the very most of the integration domain and its average is $\sim 0.363$. 

To obtain estimates for the stationary CO coverage for a given set of parameter values, we employ chains of 20 sites with periodic boundary conditions and simulate the  stochastic
 process using the lattice kinetic Monte Carlo (kMC) code {\sc kmos}\cite{hoffmann2014kmos}. We want to address stationary operation and therefore perform $10^7$ initial kMC steps for relaxation\cite{DistributionGSA}. Afterwards, we perform additional $10^7$ kMC steps and obtain an estimate of the stationary expected value by time averaging over these. We run multiple of these $2\times 10^{7}$ long trajectories with different random seeds to achieve the targeted variance.

Unlike to the previous section, where we compared relaxation of the error during the refinement, we now want to compare final errors of MLASG and ASG for different target tolerances $tol \in [{ 1,10^{-2}}]$.
As before, we choose the initial  standard deviation $\sigma_0=tol$ for MLASG. Unlike before, we employ the same small standard deviation of ${10^{-2}}$
for ASG, irrespective of the employed tolerance. In this way, we mimic
the blind usage of sparse grids without a balance between model and quadrature tolerance. For comparison, we also investigate non-adaptive full sparse grids (FSG) with just a single sample per node. As outlined in section \ref{sec:Multilevel}, this strategy 
 results from the telescopic sum approach for Monte Carlo models, if one targets at a non-adaptive quadrature. In the non-adaptive setting, this allows to reduce the quadrature error compared with highly accurate Monte Carlo estimates with same overall cost, i.e. with less nodes. If a single sample estimate has a large variance, this approach should become more efficient than the MLASG approach, because we need to control the variance for every node in MLASG. If we need very large numbers of samples to achieve this for the first node, we could equally well use the simple FSG approach with the same number of nodes. Depending on the problem, this might already lead to  negligible  quadrature errors as well as sufficiently low stochastic noise on the estimate for the integral.

\begin{figure}
\begin{center} 
    \includegraphics[width=0.45\linewidth]{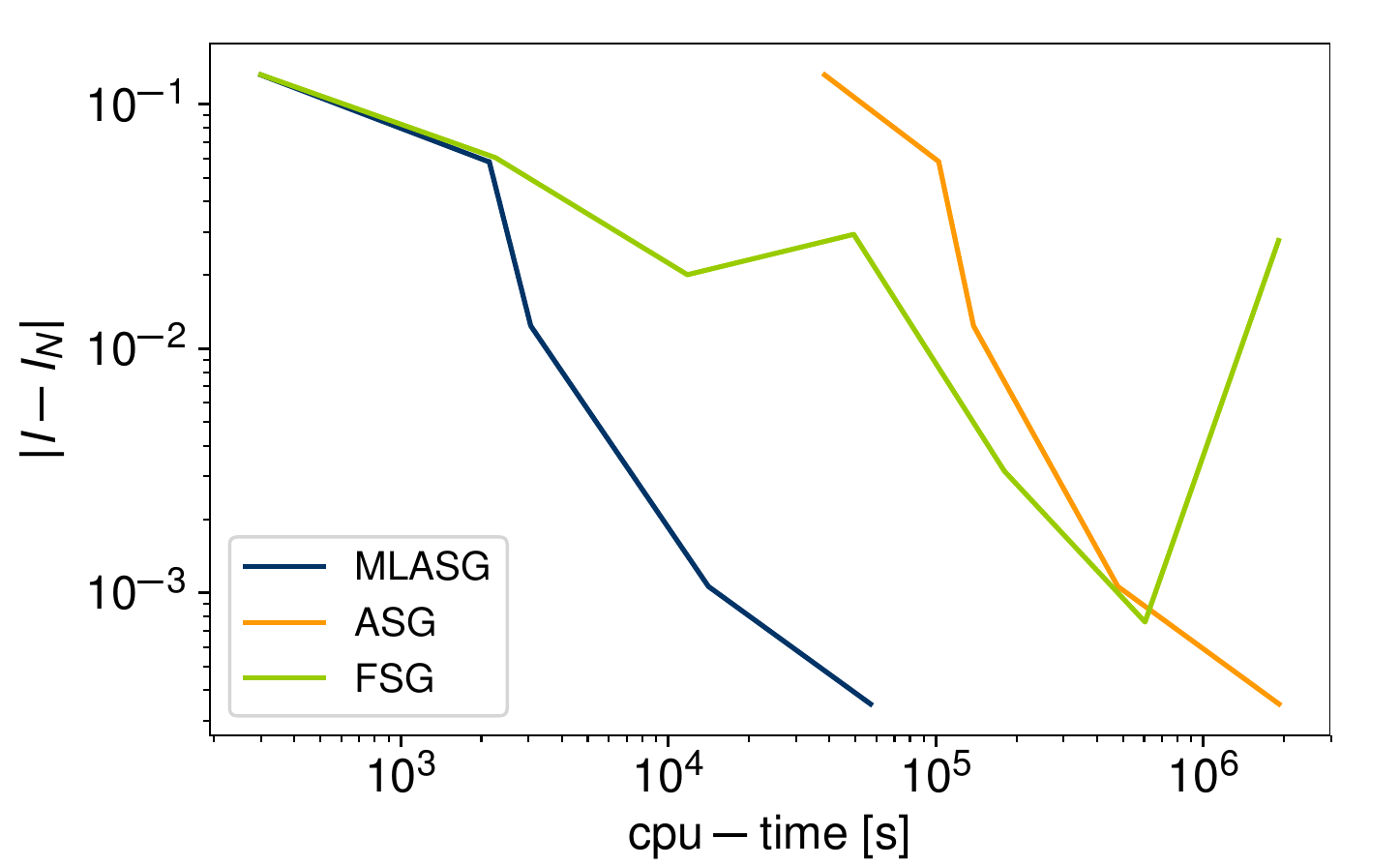}
    \includegraphics[width=0.45\linewidth]{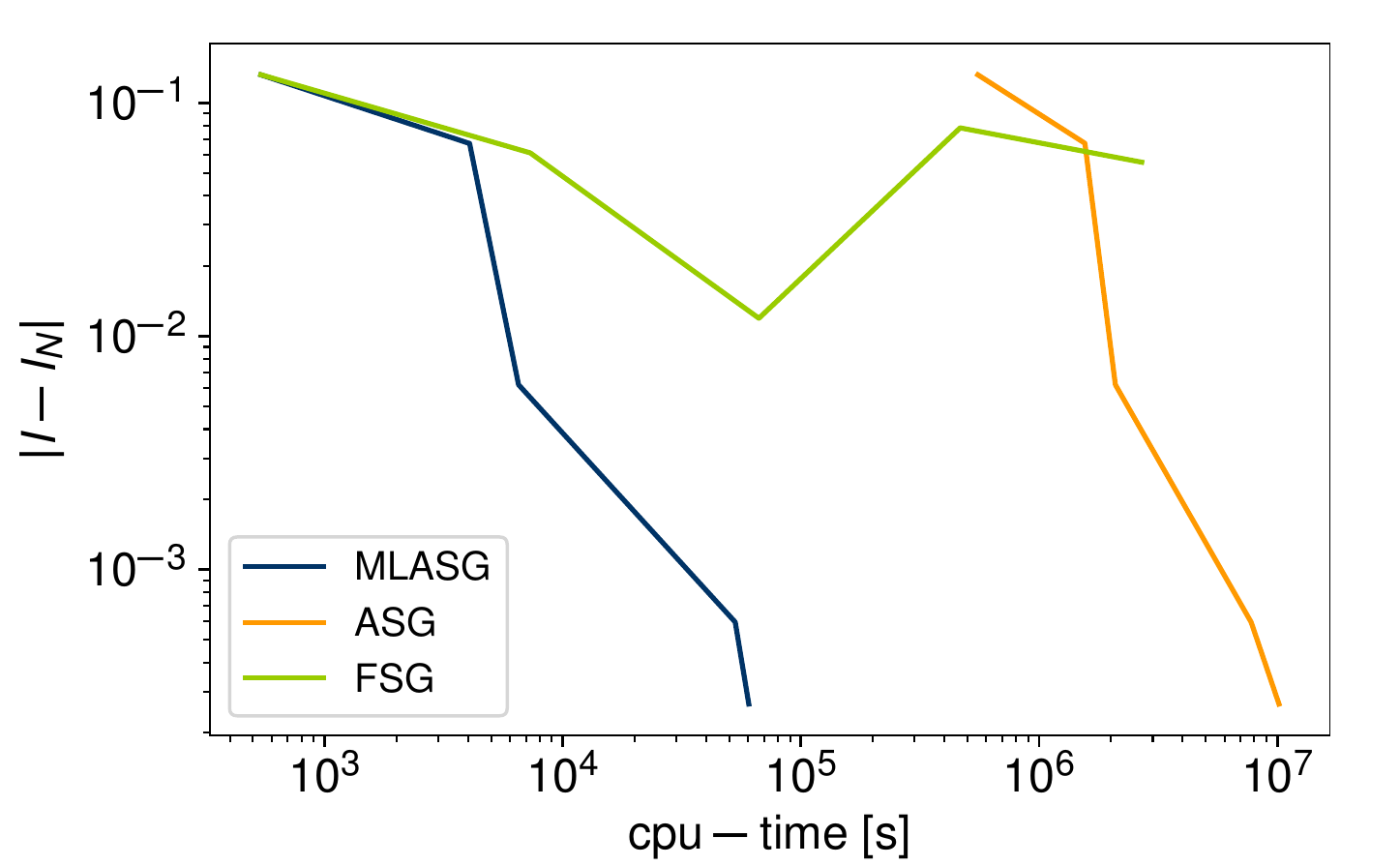}
\end{center}
\vspace{-0.5cm}
\caption{Comparison of the quadrature errors of ASG and MLASG for different tolerances $tol$ and FSG for the average CO coverage in the CO oxidation model. \textit{Left panel}: The original 7 parameter model. \textit{Right panel}: The extended 13 parameter model. }
\label{fig:ML_RuO2}
\end{figure}

The results are shown in the left panel of figure \ref{fig:ML_RuO2}, where we plot the error vs. the cpu-time spend for the kMC simulations. As expected MLASG outperforms the blind ASG and achieves a factor of $40$ speedup for the lowest tolerances - where  $\sigma_0$
is identical in both cases. ASG
seems to have a higher convergence rate. This is, however, due to the blind use of ASG. The function evaluations for low tolerances are simply unnecessarily accurate. The single sample FSG generally converges very slowly and, disregarding the incidentally accurate outliers, never reaches high accuracies. MLASG is more accurate for all computational costs and the benefit from the adaptivity therefore largely overcompensates the higher cost per node at low levels.

The benefits of the MLASG become even more pronounced, if we increase the dimensionality of the problem. To  mimic
this, we add  three
 dummy species to the model, i.e. we add the states $\ce{B^i_l},~ i \in \{1,2,3\}$ to the three states $\ce{e_l}$, $\ce{CO_l}$ and $\ce{O_l}$ of the site $l$. We assume that these can adsorb onto the surface in the same fashion as CO and add the elementary processes
\begin{align*}
 	\ce{e_l &->[k^{\rm ads}_{\rm B^i}] B^i_l},\\
  	\ce{B^i_l &->[k^{\rm des}_{\rm B^i}] e_l}
\end{align*}
where again $k^{\rm des}_{\rm B^i}=k^{\rm ads}_{\rm B^i}K_{\rm B^i}^{-1}$ and $k^{\rm des}_{\rm B^i} \in [{1 , 10^5}]$
 and $K_{\rm B^i} \in [{ 1 , 10^{4}}]$.
We thus have increased the dimensionality of the parameter space to $13$. Such additional species might, for instance, be due some trace gases is the gas phase above the catalyst, which can adsorb but do no take part in the reaction.

With this modified model, we repeat the above outlined study keeping everything else as before. 
The right panel in \ref{fig:ML_RuO2} displays the error of the methods vs. the cpu time spend for kinetic Monte Carlo sampling. The accuracy and cpu-time of the MLASG is hardly affected, which is, of course, also due to the chosen extension of the model. For the highest accuracies, MLASG is now roughly a factor $200$  more
efficient that ASG. FSG seems to hardly converge, which is expected for problems of the given dimensionality.

\section{Conclusion}
We have presented a Multilevel Adaptive Sparse Grid (MLASG) approach for the parametric numerical quadrature of Monte Carlo models. Unlike the common multilevel quadrature approaches this is not based on a rewriting of the integral as a telescoping sum. Rather, we employ the intrinsic multilevel structure of the sparse grids and control the sampling accuracy for every grid point such that the adaptive refinement strategy is not much affected by the noise in the approximate function values. For the employed basis and refinement strategy, we found that the sampling variance can be doubled with every refinement level. 
We have tested the approach on a number of examples and have found that the MLASG can achieve the same accuracy as single level Adaptive Sparse Grids (ASG) but at up to orders of magnitude lower computational costs.

In this study, we focused on quadrature, but also found hints that the higher noise level in MLASG has only a small impact on the interpolation error compared to an ASG without noise. Future research will elaborate on this finding. Also, we plan to incorporate integer programming approaches for determining
the optimal distribution of computational resources for sampling. This will allow us to test different refinement strategies but also different basis functions.

From the application point of view, we plan to test the approach on real life problems from the field of molecular simulations. 

\begin{acknowledgments}
This research was carried out in the framework of {\sc Matheon} supported by Einstein Foundation Berlin through ECMath within the sub-project SE23.
\end{acknowledgments}

\bibliography{MLASG}

\end{document}